\RequirePackage{fix-cm}
\documentclass[a4paper]{amsart}       % onecolumn (second format)
\usepackage{graphicx}

\usepackage[section]{algorithm}
\usepackage{algorithmicx}
\usepackage{algpseudocode}% http://ctan.org/pkg/algorithmicx
\usepackage{latexsym}
\usepackage{amsmath, amssymb}
\usepackage{graphicx, pgfplots}

\usepackage[colorlinks=true,linkcolor=blue,citecolor=blue,pdfpagelabels=true]{hyperref}

\usepackage{siunitx} % Table formatting
\usepackage{multirow}

\newcommand{\T}{^\mathsf{T}}
\newcommand{\cT}{\mathcal{T}}
\newcommand{\cL}{\mathcal{L}}
\newcommand{\cS}{\mathcal{S}}
\newcommand{\R}{\mathbb{R}}
\newcommand{\Ri}{R^{-1}}
\newcommand{\Mi}{M^{-1}}
\newcommand{\Mit}{M^{-\mathsf{T}}}
\newcommand{\norm}[1]{\lVert #1 \rVert}
\newcommand{\dt}{\mathrm{d}t}
\renewcommand{\exp}[1]{\mathrm{e}^{#1}}
\newcommand{\ackname}{Acknowledgements}

\begin{document}

\title[GPU acceleration of splitting schemes applied to DMEs]{GPU acceleration of splitting schemes applied to differential matrix equations}

\author[H.~Mena]{Hermann Mena}
\address{Universidad Yachay Tech, Hacienda San Jos\'{e} s/n, San Miguel de Urcuqu\'{i},  Ecuador} 
%		Universit\"{a}t Innsbruck, Technikerstra{\ss}e 13, A-6020 Innsbruck, Austria}
\email{mena@yachaytech.edu.ec}

\author[L.-M.~Pfurtscheller]{Lena-Maria Pfurtscheller}
\address{Universit\"{a}t Innsbruck, Technikerstra{\ss}e 13, A-6020 Innsbruck, Austria }
\email{Lena-Maria.Pfurtscheller@uibk.ac.at}

\author[T.~Stillfjord]{Tony Stillfjord}
\address{Max Planck Institute for Dynamics of Complex Technical Systems, Sandtorstr.\ 1, DE-39106 Magdeburg, Germany}
\email{stillfjord@mpi-magdeburg.mpg.de}

\thanks{\textit{Date:} October, 18, 2018.\\
	\indent
	2010 \textit{Mathematics Subject Classification} 65F30, 65Y05, 65F60. \\
	\indent \textit{Keywords:}{Differential {L}yapunov equations, Differential {R}iccati equations, Large-scale, Splitting schemes, GPU acceleration}}

\maketitle

\begin{abstract}
	We consider differential Lyapunov and Riccati equations, and generalized versions thereof. Such equations arise in many different areas and are especially important within the field of optimal control. In order to approximate their solution, one may use several different kinds of numerical methods. Of these, splitting schemes are often a very competitive choice. In this article, we investigate the use of graphical processing units (GPUs) to parallelize such schemes and thereby further increase their effectiveness. According to our numerical experiments, large speed-ups are often observed for sufficiently large matrices. We also provide a comparison between different splitting strategies, demonstrating that splitting the equations into a moderate number of subproblems is generally optimal.
%	\keywords{Differential {L}yapunov equations \and Differential {R}iccati equations \and Large-scale \and Splitting schemes \and GPU acceleration}
%	\subclass{ 65F30 \and 65Y05 \and 65F60}
\end{abstract}

\section{Introduction}

We are interested in differential matrix equations of Lyapunov or Riccati type, or generalized versions of these. They are all of the form
\begin{equation*}
\dot{P} = A\T P + PA + Q + G(P),
\end{equation*}
where $A \in \R^{n \times n}$ and $Q \in \R^{n \times n}$ are given matrices, $G$ is a matrix-valued function of the solution $P \in \R^{n\times n}$. For differential Lyapunov equations (DLE) we have $G(P) = 0$ and for differential Riccati equations (DRE) we have $G(P) = - P B\Ri B\T P$ with two given matrices $B \in \R^{n \times m}$ and $R \in \R^{m \times m}$. Such equations occur frequently in many different areas, such as in optimal/robust control, optimal filtering, spectral factorizations, $\mathbf{H}_\infty$-control, differential games, etc.~\cite{AboFIJ03,BasB95,IchK99,PetUS00}.

Perhaps the most relevant setting is the linear quadratic regulator (LQR) problem. There, the aim is to optimize a finite-time cost function of the form

\begin{equation*}
J(u) = \int_{0}^{T}{ \norm{y(t)}^2 + \norm{u(t)}^2 \, \dt}, \quad T \geq 0,
\end{equation*}
under the constraints that $\dot{x} = Ax + Bu$ (state equation) and $y = Cx$ (output equation, with $C \in \R^{p \times n}$). In this case, the solution to the DLE with $Q = C\T C$ gives the observability Gramian of the system, which characterizes the relevant states $x$ for the input-output mapping $u \mapsto y$. The solution of the DRE, on the other hand, provides the optimal input that minimizes $J$, in state feedback form. In fact, if $P$ solves the DRE with $Q = C^TC$ then the optimal input $u^{\text{opt}}$ is given by $u^{\text{opt}}(t) = -R^{-1}B\T P(T-t)x(t)$.

For the generalized DLE and DRE versions, an additional linear term $S P S\T$ appears in  $G(P)$, where $S \in \R^{n \times n}$ is a given matrix. Such equations also arise in the LQR setting, when a stochastic perturbation of multiplicative type is included in the state equation.

In recent years, a number of numerical methods have been suggested for large-scale DLEs, DREs and related equations. The classic ones, low-rank versions of BDF and Rosenbrock schemes~\cite{BennerMena2013,BennerMena2016,LanMS15} are usually outperformed by more modern methods such as Krylov-based projection schemes~\cite{KoskelaMena2017}, peer methods~\cite{Lang2017} or splitting schemes~\cite{Mena_etal2018,Sti15,Sti17}. In this paper, we focus on splitting schemes. These methods lower the computational cost by dividing the problem into simpler subproblems such as $\dot{P} = A\T P + PA$ and $\dot{P} = Q$ and then solve these separately, in sequence.  While the splitting of course introduces an additional error, this is generally compensated for by the decreased computational cost and leads to large speed-ups.

The hypothesis to be investigated in this paper is that utilizing a graphical processing unit (GPU) to parallelize the schemes may further greatly increase the efficiency. Such speed-ups have already been observed for other related methods for DREs~\cite{BenDEetal17,BenEMetal11,BenEMetal13} as well as for their steady-state versions: the algebraic Lyapunov and Riccati equations~\cite{BenEMetal13,BenEQetal12}. In the just mentioned cases, the basic building block of the schemes is the computation of the matrix sign function, which requires the inversion of a large dense matrix. In a splitting scheme, the basic building block is instead the computation of the action of a matrix exponential on a skinny matrix. Speed-ups have previously been observed for applications where matrix exponentials are multiplied by vectors~\cite{Auer_etal2018,EinkemmerOstermann2013}, see also~\cite{Farquhar_etal2016}. In these works, a speed-up is generally not observed for ``small'' matrices ($n \lesssim 1000$), and the speed-up is of limited size when the matrices are sparse rather than dense. As we are typically interested in at least medium-sized problems ($1000 \lesssim n \lesssim 10000$) we do expect to see a significant speed-up. Moreover, while we are necessarily considering the sparse case, we are not simply computing the action of the matrix exponential on vectors, but on skinny block matrices. This increases the parallelizability of the problem and makes the sparsity issues noted in e.g.~\cite{BellGarland2008,EinkemmerOstermann2013,Goumas_etal_2008} less relevant.

Since the relevant methods are mainly implemented in Matlab, we restrict ourselves to utilizing its built-in GPU support~\cite{ReeZ12} via NVIDIA's CUDA~\cite{Nickolls_etal2008} parallel programming interface. We do not claim that this approach leads to the best possible performance. The point is rather to demonstrate that quite simple changes to the implementations of the splitting schemes may lead to much better performance, when one has access to a GPU. Our results already show a remarkable improvement in efficiency, and this can only increase with further optimisations and the use of more advanced techniques tailored to specific problems.

In addition, we provide comparisons between different splitting strategies for DLEs and DREs. We particularly address questions that naturally arise while solving these equations by splitting methods. E.g., should the DLE be split at all? Should the DRE be split into two or three subproblems? Our results in this direction demonstrate that it is usually beneficial to use the smallest number of splits. However, when $Q$ is sufficiently small it is beneficial to split it  too, since the extra error is similarly small and the subproblems $\dot{P} = A\T P + PA$ and $\dot{P} = Q$ are very cheap to compute compared to $\dot{P} = A\T P + PA + Q$.

An outline of the article is as follows. In Section~\ref{sec:splitting_schemes} we review the idea behind splitting schemes and apply them to all the mentioned equation types. Then we consider implementation details in Section~\ref{sec:implementation}. The simple changes necessary for GPU utilization, and a discussion on what efficiency improvements may be expected is given in Section~\ref{sec:GPU}.
The actual speed-ups are presented in Section~\ref{sec:experiments}, in the form of several numerical experiments. Finally, we summarize our conclusions in Section~\ref{sec:conclusions}.

\section{Splitting schemes} \label{sec:splitting_schemes}
Splitting schemes are numerical methods that are applicable to differential equations that have a natural decomposition into two (or more) parts;
\begin{equation*}
  \dot{P} = F(P) = F_1(P) + F_2(P), \quad P(0) = P_0.
\end{equation*}
With ``natural decomposition''  we mean that the subproblems
\begin{equation*}
  \dot{P} = F_1(P), \quad \text{and} \quad  \dot{P} = F_2(P)
\end{equation*}
are either simpler or cheaper to solve than the full problem $\dot{P} = F(P)$. This is the case in many problems, with the most common example being reaction-diffusion equations $\dot{x} = \Delta x + f(x)$ with homogeneous Dirichlet boundary conditions. In this case, there are highly optimized methods for the pure diffusion problem $\dot{x} = \Delta x$, while the subproblem $\dot{x} = f(x)$ often turns into a local rather than global problem --- i.e.\ it is enough to solve $\dot{x_i} = f(x_i)$ for every discretization point $x_i$. (For some caveats in the case of other boundary condition types, see e.g.~\cite{AloCR17,EinO15,EinO16}.)
In the following, we denote the solution to $\dot{P} = F_k(P)$, $P(0) = P_0$, by $P(t) =: \cT_k(t)P_0$.

The most basic and commonly used (exponential) splitting schemes are the Lie and Strang splittings. They are given by the time stepping operators
\begin{equation*}
   \cL_h P_0 =  \cT_2\left(h\right) \cT_1\left(h\right) P_0, \quad \text{and} \quad \cS_h P_0 = \cT_1\left(\frac{h}{2}\right) \cT_2\left(h\right) \cT_1 \left(\frac{h}{2}\right) P_0,
\end{equation*}
respectively, where $h$ is the time step. Of course, the roles of $\cT_1$ and $\cT_2$ might be interchanged. The schemes are then defined by 
\begin{equation*}
  P^L_{k+1} = \cL_h P^L_k, \quad \text{and} \quad P^S_{k+1} = \cS_h P^S_k, 
\end{equation*}
with  $P^L_0 = P^S_0 = P_0$. Here, $P^L_k$ and $P^S_k$ both approximate $P(kh)$. The Lie splitting is first-order accurate while Strang splitting is  second-order accurate under certain conditions on $F_1$, $F_2$ and $F$, see e.g.~\cite{HunV03}. For simplicity, we restrict ourselves to the Strang splitting scheme in this paper, but one might also consider higher-order schemes~\cite{DeRS16,HanO09,Sti17}, or schemes where the subproblems are not solved exactly, see e.g.~\cite{HanS14,HunV03}.

Clearly, one might continue the splitting procedure if the system is naturally decomposed into more than two parts. If 
\begin{equation*}
  \dot{P} = F_1(P) + F_2(P) + F_3(P),  
\end{equation*}
then applying the Lie and Strang splitting schemes twice leads to the schemes
\begin{align*}
   \tilde{\cL}_h P_0 &=  \cT_3\left(h\right)  \cT_2\left(h\right) \cT_1\left(h\right) P_0, \quad \text{and} \quad \\
\tilde{\cS}_h P_0 &= \cT_1\left(\frac{h}{2}\right)  \cT_2\left(\frac{h}{2}\right)  \cT_3\left(h\right)  \cT_2\left(\frac{h}{2}\right) \cT_1 \left(\frac{h}{2}\right) P_0.
\end{align*}
Again, the roles of $\cT_1$, $\cT_2$ and $\cT_3$ might be interchanged. Different compositions with a possibly higher number of operators might also be considered, in order to optimize the structure of the error. We refer to~\cite{AuzKT15} but do not consider such methods here.

Like essentially every other method for solving differential matrix equations, the splitting schemes need to make use of low-rank structure in order to be competitive in the large-scale setting. This means that we can expect the singular values of the symmetric, positive semi-definite solution $P$ to decay rapidly, see e.g.~\cite{AntSZ02,BakES14,BenB13,Pen00,SorZ02}, and thus we can factorize $P \approx L D L\T$ for $L \in \R^{n \times r}$, $D \in \R^{r \times r}$ with $r \ll n$. By formulating the methods to only operate on $L$ and $D$ and never explicitly form the product $LDL^T$, we drastically lower both the memory requirements and the computational cost.

In the following, we outline different splitting strategies for all the matrix equations mentioned so far, and also review how to low-rank-factorize each arising subproblem.

\subsection{Differential Lyapunov equations} \label{subsec:DLE} 
As a first example, we consider the differential Lyapunov equation
\begin{align}
\dot{P} = A\T P + PA + Q, \quad P(0) = P_0   \label{eq:lyap}
\end{align} 
Here we may choose $F_1$ as the linear part and $F_2$ as the constant term, i.e.
\begin{equation*}
  F_1(P) = A\T P + P A, \quad \text{and} \quad F_2(P) = Q .
\end{equation*}

These subproblems can be solved explicitly and the solutions at time $h$ are given by
\begin{align*}
  \cT_1(h) P_0 &= \exp{hA\T} P_0 \exp{hA},\\
  \cT_2(h) P_0 &= P_0 + h Q . 
\end{align*}
It is easily seen that if we have the $LDL\T$-factorizations $P_0 = LDL\T$ and $Q = L_QD_QL\T_Q$, then we can also factorize these solutions as
\begin{align}
  \cT_1(h) P_0 &= \left(\exp{hA\T} L\right) D \left(\exp{hA\T} L \right)\T, \label{eq:F_sol_LDL}\\
  \cT_2(h) P_0 &= \begin{bmatrix} 
                    L & L_Q
                 \end{bmatrix}
                 \begin{bmatrix}
                     D & 0   \\
                     0 & hD_Q
                 \end{bmatrix}
                 \begin{bmatrix} 
                    L & L_Q
                 \end{bmatrix}\T. 
\end{align}
We could also note that the exact solution to the full problem is given by
\begin{equation}
  P(t) = \exp{tA\T} P_0 \exp{tA} + \int_{0}^{t}\exp{sA\T} Q \exp{sA} ds,  \quad t \in [0,T], \label{eq:full}
\end{equation} 
where the integral term may be approximated by high-order quadrature as in~\cite{Sti15}. While this does not result in a splitting scheme of the form described above, we still include it in our experiments due to its similarity and efficiency.

\subsection{Differential Riccati equations}\label{subsec:DRE}
A second example is given by the differential Riccati equation:
\begin{align}
\dot{P} = A\T P + PA + Q - P B\Ri B\T P, \quad P(0) = P_0 . \label{eq:ricc}
\end{align} 
In this case, we can either split in three terms;
\begin{equation*}
  F_1(P) = A\T P + P A, \quad F_2(P) = Q, \quad \text{and} \quad  F_3(P) = -P B \Ri B\T P,
\end{equation*}
or two terms\footnote{We deliberately use $F_{12}$ and $F_3$ here rather than $F_1$ and $F_2$, in order to not change the meaning of the previously defined $F_1$ and $F_2$. The two-term splitting schemes are obviously still well-defined after substituting the proper numbers.
};
\begin{equation*}
  F_{12}(P) = A\T P + P A + Q, \quad \text{and} \quad  F_3(P) = -P B\Ri B\T P.
\end{equation*}
The latter was advocated in~\cite{Sti15,Sti17} because (experimentally) the error constant in the three-term splitting is much larger. However, the three-term splitting does not need to approximate the integral term, and thus the larger error might be compensated for by a lower computational cost.

In either case, we note that the solution  at time $h$ to the problem ${\dot{P} = F_3(P)}$, $P(0) = P_0$, is given explicitly by
\begin{equation}\label{eq:nonlinear}
\cT_3(h) P_0 = (I + h  P_0 B\Ri B\T)^{-1} P_0. 
\end{equation}
A low-rank factorization is given by
\begin{equation*}
  \cT_3(h) LDL\T = L (I + h DL\T B \Ri B\T L)^{-1}DL\T.
\end{equation*}
Note that the $I$ in this equation is not the same identity matrix as in the previous equation, because the $L$-part of $P_0$ has moved. We thus only need to solve a small linear equation system.

\subsection{Generalized Lyapunov equations} \label{subsec:gen_DLE}
We further consider a generalized Lyapunov equation of the form
\begin{align}
  \dot{P} = A\T P +PA + Q + S P S\T, \quad P(0) = P_0. \label{eq:lyapGen}
\end{align} 
We again split the equation and obtain three subproblems defined by\footnote{For the same reason as in the previous note, we use $F_4$ rather than $F_3$ here.} 
\begin{equation*}
  F_1(P) = A\T P + P A, \quad F_2(P) = Q, \quad \text{and} \quad  F_4(P) = S P S\T.
\end{equation*}
The first two subproblems are handled as before, whereas we approximate $\cT_4(h)(P)$ by the midpoint rule, analogously to what is done in~\cite{DamMS17}:
\begin{align*}
  \cT_4(h)P_0 \approx P_0 + h S \Big(P_0 + \frac{h}{2} S P_0 S\T\Big) S\T .
\end{align*} 
Given $P_0 = L_0 D_0 L_0\T$, we get $\cT_4(h)P_0 \approx LDL\T$, where 
\begin{equation*}
 L = \left[L_0, \sqrt{h} SL_0, \frac{h}{\sqrt{2}} S^2 L\right] , \quad \text{and} \quad D = \texttt{blkdiag}(D_0,D_0,D_0) ,
\end{equation*}
where \texttt{blkdiag} is the block diagonal operator that puts its block arguments on the diagonal of an otherwise zero matrix.

We note that when using a second-order splitting scheme like the Strang splitting, it is necessary to use a second-order method like the midpoint rule in order to preserve the overall convergence order. If we use instead a first-order scheme like the Lie splitting, it is sufficient to approximate $\dot{P} = F_4(P)$ by e.g.\ the explicit Euler method.

\subsection{Generalized Riccati equations}\label{subsec:gen_DRE}
Moreover, we study a generalized Riccati equation given by
\begin{align}
  \dot{P} = A\T P + PA + Q + S P S\T - P B\Ri B\T P , \quad P(0) = P_0 ,
\end{align} 
and split this equation into three subproblems of the form
\begin{equation*}
  F_{12}(P) = A\T P + P A + Q , \quad F_3(P) = -P B\Ri B\T P, \quad \text{and} \quad  F_4(P) = S P S\T.
\end{equation*}
These subproblems are solved similarly as in the previous subsections. We do not consider a four-term splitting since experience suggests that the extra error due to the splitting would become prohibitively large.

\section{Implementations} \label{sec:implementation}
In this section we describe the implementation of the Strang splitting scheme applied to the differential matrix equations discussed in Section~\ref{sec:splitting_schemes}. Other splitting schemes such as the Lie splitting are implemented analogously.
 
In all the considered equations, the most demanding part is to compute the action of the matrix exponential in~\eqref{eq:F_sol_LDL} efficiently. In~\cite{CalKOR14,CalKOR16} the authors considered an algorithm based on Leja interpolation and showed that applying the algorithm to a matrix derived from a spatial discretization of a differential operator is very efficient. We therefore use this method to compute $\exp{hA}L$ for different skinny matrices $L$, and denote it by \texttt{expleja} in the following. 

First, we consider the DLE case. The discussion in Section~\ref{subsec:DLE} immediately leads to Algorithm~\ref{alg:DLE_splitting}.
\begin{algorithm}[!htb]
\caption{Solving DLE by Strang splitting}
\label{alg:DLE_splitting}
\begin{algorithmic}[1]
\State Given: $A$, $Q$, $P_0$, $T$, $N_t $, $h = \frac{T}{N_t}$.  
\State Compute $LDL\T$-decompositions of $Q = L_Q D_Q L_Q\T$ and $P_0 = LDL\T$.
\State Compute parameters \texttt{param} for Leja interpolation. 
\For{$k = 1, \ldots, N_t$}
\State $L = \texttt{expleja}(h/2,A,L,\texttt{param})$
\State $L = [L, L_Q]$
\State $D = \texttt{blkdiag}(D, hD_Q)$; 
\State $[L,D] = \texttt{column\_compression}(L,D) $;
\State $L = \texttt{expleja}(h/2,A,L,\texttt{param})$
\EndFor
\State $P = LDL\T;$
\end{algorithmic}
\end{algorithm}

On the other hand, as mentioned in Subsection~\ref{subsec:DLE} it is possible to derive an explicit form of the solution of the DLE given by~\eqref{eq:full}. Following~\cite{Sti17}, we use a high-order quadrature rule to compute an approximation to the integral term. This computation is again based on using \texttt{expleja}, now to compute $\exp{s_kA}L_Q$ for various $s_k \in [0,h]$ with the $LDL\T$-factorization $Q = L_QD_QL_Q\T$. This leads to the alternative Algorithm~\ref{alg:DLE_quad_rule}, which (as noted in Section~\ref{subsec:DLE}) is not a splitting scheme per se.

\begin{algorithm}[!htb]
\caption{Solving DLE by quadrature}
\label{alg:DLE_quad_rule}
\begin{algorithmic}[1]
\State Given: $A$, $Q$, $P_0$, $T$, $N_t $, $h = \frac{T}{N_t}$. 
\State Repeat Steps 2 and 3 from Algorithm~\ref{alg:DLE_splitting}.
\State Approximate integral: 
	\begin{itemize}
	\State Compute $n$ nodes $s_k$ and weights $w_k$ of a quadrature formula;
	\State $L_I = [\texttt{expleja}(s_1, A, L_Q), \ldots \texttt{expleja}(s_n, A, L_Q) ]$;
	\State $D_I = \texttt{blkdiag}(w_1 D_Q, \ldots, w_n D_Q)$;
	\State $[L_I,D_I] = \texttt{column\_compression}(L_I,D_I) $.
	\end{itemize} 
\For{$k = 1, \ldots, N_t$}
\State $L = [\texttt{expleja}(h,A,L,\texttt{param}),L_I]$
\State $D = \texttt{blkdiag}(D,D_I);$
\State $[L,D] = \texttt{column\_compression}(L,D) $;
\EndFor
\State $P = LDL\T$ .
\end{algorithmic}
\end{algorithm}

We note that in both Algorithm~\ref{alg:DLE_splitting} and Algorithm~\ref{alg:DLE_quad_rule} there is a so-called column compression step. This refers to the procedure of discarding (almost) linearly dependent columns from $L$, and serves to keep the number of columns in the approximations small. Without such a step, each iteration of Algorithm~\ref{alg:DLE_splitting} (for example) would add the columns in $L_Q$ to $L$, while the rank would likely stay similar. The compression can be performed in various ways, usually by computing either a reduced rank-revealing QR factorization or a reduced SVD~\cite{LanMS15}. 
Here, we employ a reduced SVD factorization, followed by a diagonalization of the small resulting system. It is cheap as long as the rank of the solution stays low, which is the case in many applications.

As noted in Section~\ref{sec:splitting_schemes}, we also want to approximate the solutions to DREs and generalized DLEs and DREs. Therefore, we further have to solve the subproblems given by $F_3$ and $F_4$. Pseudo-codes for these computations, based on the low-rank factorizations given in Sections~\ref{subsec:DRE}--\ref{subsec:gen_DLE}, are shown in Algorithms~\ref{alg:Ricc_term}~--~\ref{alg:bil_term}. 

\begin{algorithm}[!htb]
\caption{Solving $\dot{P} = F_3(P)$ over $[0, h]$}
\label{alg:Ricc_term}
\begin{algorithmic}[1]
\State Given: $B$, $\Ri$, $h$ and a low-rank factorization of $P = LDL\T$.
\State Compute $D = (I + hDL\T B\Ri L)^{-1} D$;
\State $P = LDL\T$. 
\end{algorithmic}
\end{algorithm}

\begin{algorithm}[!htb]
\caption{Solving $\dot{P} = F_4(P)$ over $[0, h]$}
\label{alg:bil_term}
\begin{algorithmic}[1]
\State Given: $S$, $h$ and a low-rank factorization of $P = LDL\T$.
\State Compute $L = [L,\sqrt{h} S L, h / \sqrt{2} \, S^2 L  ];$
\State Compute $D = \texttt{blkdiag}(D,D,D)$;
\State $[L,D] = \texttt{column\_compression}(L,D) $;
\State $P = LDL\T$.  
\end{algorithmic}
\end{algorithm}

We use three approaches to split the DRE: First, we incorporate Algorithm~\ref{alg:Ricc_term} in Algorithm~\ref{alg:DLE_quad_rule} in order to solve the Lyapunov part of the equation via quadrature and the nonlinear term via the exact solution formula in Equation~\ref{eq:nonlinear}, forming 
\begin{equation*}
	\cT_{12}\left(\frac{h}{2}\right) \cT_3\left(h\right) \cT_{12} \left(\frac{h}{2}\right) P_0 .
\end{equation*}
 Further, we consider the three-term splitting  
\begin{equation*}
	\cT_1\left(\frac{h}{2}\right) \cT_2\left(\frac{h}{2}\right)  \cT_3\left(h\right)  \cT_2\left(\frac{h}{2}\right) \cT_1 \left(\frac{h}{2}\right) P_0 ,
\end{equation*}
by extending Algorithm~\ref{alg:DLE_splitting} with a third step given by Algorithm~\ref{alg:Ricc_term}.
 Finally we reverse the order of the three-term splitting 
\begin{equation*}
	\cT_1\left(\frac{h}{2}\right) \cT_3\left(\frac{h}{2}\right)  \cT_2\left(h\right)  \cT_3\left(\frac{h}{2}\right) \cT_1 \left(\frac{h}{2}\right) P_0 .
\end{equation*}
Due to the additional splitting term, further errors are introduced, but since the integral does not have to be approximated the three-term splitting codes are less computationally demanding than the two-term splittings. 

The generalized DLE can be solved by the same three approaches. Using Algorithm~\ref{alg:bil_term}, $\cT_3$ is replaced by $\cT_4$ in the previous three formulas.
Finally, we consider a three-term Strang splitting for the generalized DRE, given by
\begin{equation*}
	\cT_{12}\left(\frac{h}{2}\right) \cT_3\left(\frac{h}{2}\right)  \cT_4\left(h\right)  \cT_3\left(\frac{h}{2}\right) \cT_{12} \left(\frac{h}{2}\right) P_0 .
\end{equation*}  

The modifications to Algorithms~\ref{alg:DLE_splitting} and~\ref{alg:DLE_quad_rule} for the (generalized) DRE and generalized DLE cases through use of Algorithms~\ref{alg:Ricc_term} and~\ref{alg:bil_term} are obvious, and we therefore omit full listings of these versions.

\section{GPU considerations} \label{sec:GPU}
All the algorithms in the previous section were implemented in Matlab. For GPU acceleration we used the Parallel Computing Toolbox, which interfaces with the CUDA library. This is a framework for general purpose computing on GPUs. Recent releases of Matlab expose a large fraction of this framework as overloaded built-in functions, i.e. precompiled code that operates either on the CPU or the GPU, depending on where the data currently resides. Thus, e.g., solving a system of linear equations $Ax = b$ on the GPU can be accomplished by the familiar syntax \texttt{x = A\textbackslash b} after $A$ and $b$ have been instantiated as objects on the GPU. This data transfer is performed by the \texttt{gpuArray} function. The result $x$ may then be transferred back to the CPU by use of the \texttt{gather} function. We refer to e.g., \cite{ReeZ12}.
In general, communication between the CPU and the GPU is expensive. We therefore first move all the data to the GPU, do all vector- and matrix-computations on the GPU using built-in functions and transfer only scalar quantities and the final results back to the CPU.

The main computational effort in all the algorithms is the computation of the matrix exponential actions via the \texttt{expleja} code. This consists of a (one-time) estimation of the spectrum of $A$ by the Gershgorin disk theorem, a (one-time) computation of exponential interpolation parameters, and a Newton interpolation~\cite{CalKOR14,CalKOR16}. These functions all depend only on matrix-vector or matrix-matrix products and simple built-in functions like \texttt{diag}, \texttt{sum} and \texttt{abs}, all of which have overloaded GPU-versions. There is thus no need for any changes to the main code.

\subsection{Main routines and limiting factors} \label{subsec:main_routines}
As will be demonstrated in Section~\ref{sec:experiments}, around 90-95\% of the total computation time is spent in the Newton interpolation part of the \texttt{expleja} code. On the CPU side, this can be further broken down into the multiplication of a sparse matrix with a dense skinny matrix (65-75\%), the computation of the 1-norm of a dense skinny matrix (15-25\%) and the addition of two skinny matrices (5\%). On the GPU side, the ranking of these operations are typically the same, but the relative percentages differ.

All of these operations are memory-bound, i.e., their computation is limited by memory bandwidth rather than processing power. This can easily be confirmed by considering the number of necessary read/writes compared to the number of actual floating point operations.

\subsection{Possible perfomance gains} \label{subsec:possible_gains}
 Since modern GPUs feature larger memory bandwidths than comparable CPUs and since the main operations are also highly parallelizable, we expect to see a speed-up when utilizing the GPU. This speed-up will likely not be as large as for a compute-bound problem, where GPUs excel, but should still be significant, especially considering the essentially zero cost of extra implementation effort.

If both the GPU and CPU operated at peak performance, the observed speed-up would simply be the ratio of the respective memory bandwidths. This will, however, not be the case in practice. Still, one might expect that both platforms operate at a similar percentage of peak performance, and that the ratio will stay similar. In practice, however, this will also not be the case, due to differences in code optimization. In the current application, the overwhelming majority of the computations are performed in very low-level operations which we can not influence. Since Matlab is not open source, we have no insight into what particular algorithm is used or how well it is optimized. 
Because the main operations are memory-bound, efficient memory allocation also plays a large role. Here, again, we have no insight into what strategies Matlab follows. Finally, the CPU typically also has one additional layer of (larger) cache memory than the GPU, which further complicates things.
For these reasons, it is difficult to predict what kind of speed-up to expect.

\section{Numerical experiments} \label{sec:experiments}
The aim of this section is to apply the different splitting strategies to various examples and demonstrate that the GPU implementation consistently outperforms the CPU version, often by a large margin.

We first describe a number of examples, including two arising from real-world problems. As an implementation verification, we then test our codes on the first, small-scale, problem, where we can compute an accurate reference solution by vectorization of the problem. We observe the correct orders of convergence and also verify that the GPU and CPU codes indeed give the same results. Then, we compare the speed of the two platforms by applying the different algorithms to the given test examples, and demonstrate that GPU acceleration is advantageous in all cases.

The tests were run on two different systems. The first one, hereafter referred to as ``System 1'', has an Intel Xeon E5-2630v3 CPU and a Tesla K80. The K80 contains two separate GPUs, of which we use only one. The maximum memory bandwidths are here 59 GB/s\footnote{\url{https://ark.intel.com/products/83356/Intel-Xeon-Processor-E5-2630-v3-20M-Cache-2_40-GHz}} for the CPU and 240 GB/s\footnote{\url{https://www.nvidia.com/en-us/data-center/tesla-k80/}} for the GPU. This system has $24$ GB RAM available.
The second system, hereafter referred to as ``System 2'', is one node of the Mechthild\footnote{\url{http://www.mpi-magdeburg.mpg.de/cluster/mechthild}} HPC cluster at the Max Planck Institute Magdeburg. This has an Intel Xeon Silver 4110 CPU and a Tesla P100 GPU. The maximum memory bandwidths are here 115 GB/s\footnote{\url{https://ark.intel.com/products/123547/Intel-Xeon-Silver-4110-Processor-11M-Cache-2_10-GHz}} for the CPU and 732 GB/s\footnote{\url{https://www.nvidia.com/en-us/data-center/tesla-p100/}} for the GPU. This system has $192$ GB RAM available.

In all our experiments, we use the tolerance $10^{-16}$ for both the column compression and the Leja interpolation. This ensures that the approximations are not unnecessarily truncated, and that the matrix exponential actions are essentially exact. We use Matlab R2017a on System~1 and R2017b on System~2. In both cases, we deactivate the Java Virtual Machine by \texttt{-nojvm}. The computing times of the CPU algorithms are estimated by the command \texttt{tic - toc}. For the GPU algorithms we do the same, except that we also call the \texttt{wait} function to ensure that all threads on the GPU have finished their computations before the measurements.

\subsection{Experiment descriptions}\label{subsec:description_experiments}

\subsubsection[Example 1]{Example 1: Heat equation random model}\label{example:one}
We first consider the Laplacian on the unit square with homogeneous Dirichlet boundary conditions. By discretizing it using central second-order finite differences with $n_x$ grid points in each space dimension, we acquire a matrix $A \in \R^{n\times n}$ with $n = n_x^2$. We let $Q$ and $P_0$ be randomly chosen matrices of rank $2$ and rank $5$, respectively and take $B$ to be a randomly chosen matrix of size $n \times 1$. This corresponds to optimal (distributed) control of the heat equation, with a single input and two outputs, and gives rise to a DRE. By ignoring the $B$ matrix, we get instead a DLE where the solution corresponds to the time-limited Gramian of the system. In both cases, we use the final time $T = \frac{1}{2}$.

\subsubsection[Example 2]{Example 2: Stochastic heat transfer}\label{example:two}
For the generalized matrix equations, we consider an example introduced in \cite{BenD11} arising from a stochastic heat transfer problem. The matrix $A$ again denotes the discretized $2D$ Laplacian on the unit square, but now with homogeneous Dirichlet boundary conditions on two edges. On the third edge, we implement control through the fixed boundary condition $x = u$, and on the final edge a stochastic Robin boundary condition $n \cdot \nabla x = 0.5(0.5 + dW)x$ is applied, where $W$ is a Brownian motion. This leads to a matrix $B \in \mathbb{R}^{n\times 1}$ and a (sparse) matrix $S \in \mathbb{R}^{n\times n}$. The matrix  $Q = CC\T$ is defined by letting $C = \frac{1}{n} (1,\ldots, 1)$ be the matrix representation of the mean. Similarly to the previous example, we may acquire a generalized DLE instead by simply ignoring the matrix $B$. (Then there is a homogeneous Dirichlet boundary condition also on the third edge.) In both cases, we use the final time $T = \frac{1}{2}$.

\subsubsection[Example 3]{Example 3: Simulation of El Ni\~{n}o}\label{example:three}
As a third example, we consider the real-world weather phenomenon El Ni\~no. This is characterized by an unusual warming of the sea surface temperature in the Indo-Pacific ocean. It can be modeled by a stochastic advection equation driven by additive noise~\cite{PenSar95} and its covariance is given by a DLE of the form
\begin{equation*}
  \dot{P}(t) = AP(t) + P(t) A\T + Q,
\end{equation*}
see~\cite{Mena_etal2018,MenaPf2017} for details. The matrix $A$ arises from a centered finite difference approximation of the advection operator and $Q$ is the discretized covariance operator of the random noise. We consider here the different discretization resolutions corresponding to $n = 624,3900,7800$ and $15600$ and use the final time $T = 100$. We note that this problem only yields to a DLE.

\subsubsection[Example 4]{Example 4: Simulation of steel cooling}\label{example:four}
For our final example we consider the optimal cooling of steel profiles. This problem has been widely studied in the literature, for details see \cite{BenS05,Saa03}. It gives rise to a DRE of the form
\begin{equation*}
	M\T \dot{P} M = A\T P M + M\T P A + Q - M\T P B\Ri B\T P M .
\end{equation*}
The matrices $M$ and $A$ are the mass and stiffness matrices resulting from a finite element discretization of the Laplacian on a non-convex polygonal domain (the steel profile). Q is chosen as $C\T C$, where $C$ is the discretization of an operator that measures temperature differences between different points in the domain. (We want an even temperature distribution.) Finally, the matrix $B$ is the discretization of the operator that implements the Neumann boundary conditions of the Laplacian -- this results in a boundary control application.
Cancelling $M\T$ and $M$ leads to the equation
\begin{equation*}
	\dot{P}  = \Mit A\T P + P A \Mi + \Mit Q \Mi -P B\Ri B\T P, 
\end{equation*}
which we can treat as outlined in Section~\ref{subsec:DRE} after replacing $A$ by $\Mi A$ and $Q$ by $\Mit Q \Mi$.

We note that we would normally never explicitly compute the (generally dense) matrices involving $\Mi$. In the CPU code, we form and reuse an incomplete LU decomposition of $M$ to cheaply solve a linear equation system whenever the action of $\Mi$ or $\Mit$ is required. In the GPU code, the issue is unexpectedly complicated by the fact that Matlab's CUDA interface does not support solving equation systems with sparse system matrices and (dense) block right-hand sides. This \emph{is} supported in the cuSPARSE library of CUDA itself, so until the Matlab interface is extended one might theoretically implement this capability by a MEX extension. In order to demonstrate performance gains by rather easy means, however, we do not do this. Instead, we compute and store a dense LU factorization. This is clearly not viable for truly large-scale problems, but problems of up to size $n \approx 3 \cdot 10^{4}$ are easily possible on our available hardware, and up to $n \approx 5.5 \cdot 10^{4}$ if $A \Mi$ is explicitly formed at a slightly higher initial cost. Despite the heavy additional memory requirement, the GPU parallelization will lead to a significant speed-up.

An additional issue related to the mass matrix is the original Leja point interpolation method for the computation of matrix exponential actions. One of the main steps of this algorithm computes an estimate of the spectrum of $A$ by the use of Gershgorin discs. Since computing these requires direct access to the elements in $A$, it is not directly applicable to $A \Mi$ without explicitly forming the matrix. To get around this issue and still acquire a cheap estimate, we utilized the results of~\cite{Nak11} which extends the Gershgorin approach to generalized eigenvalue problems. In our experience, this method overestimates the imaginary part of the spectrum but otherwise works well. We note that if the GPU code utilizes dense matrices, we may of course simply compute $A \Mi$ and apply the original Leja point method. Since we expect to be able to work with sparse matrices in the near future, however, we follow the approach outlined above in both the CPU and GPU codes.

In the following we consider the discretizations corresponding to $n =  371$, $1357$, $5177$ and $20209$, for which the matrices have been pre-computed . We take $\Ri = I$, $P(0) = 0$ and integrate until $T = 450$.

\subsection{Implementation verification} \label{subsec:implementation_verification}
In order to verify our implementations, we investigate the convergence properties of the methods when applied to~\nameref{example:one} and~\nameref{example:two} with $n = 25$. The reference solutions are computed by vectorizing the system and applying the Matlab routine \texttt{ode15s} with relative tolerance $2.22 \cdot 10^{-14}$ (which is the lowest possible relative tolerance) and absolute tolerance $10^{-20}$. We show only the results from the GPU versions of the code on System 1 to minimize clutter, but the CPU versions yield the same results and so do the simulations on System 2.

\begin{figure}[!htb]
  \centering
  \includegraphics[width=0.49\textwidth]{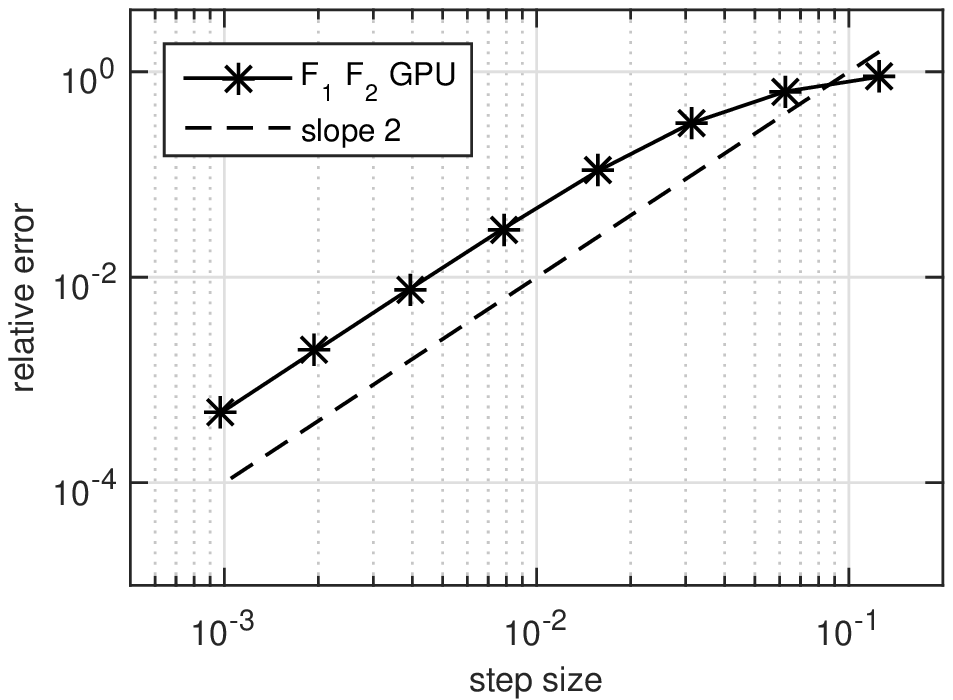}
  \includegraphics[width=0.49\textwidth]{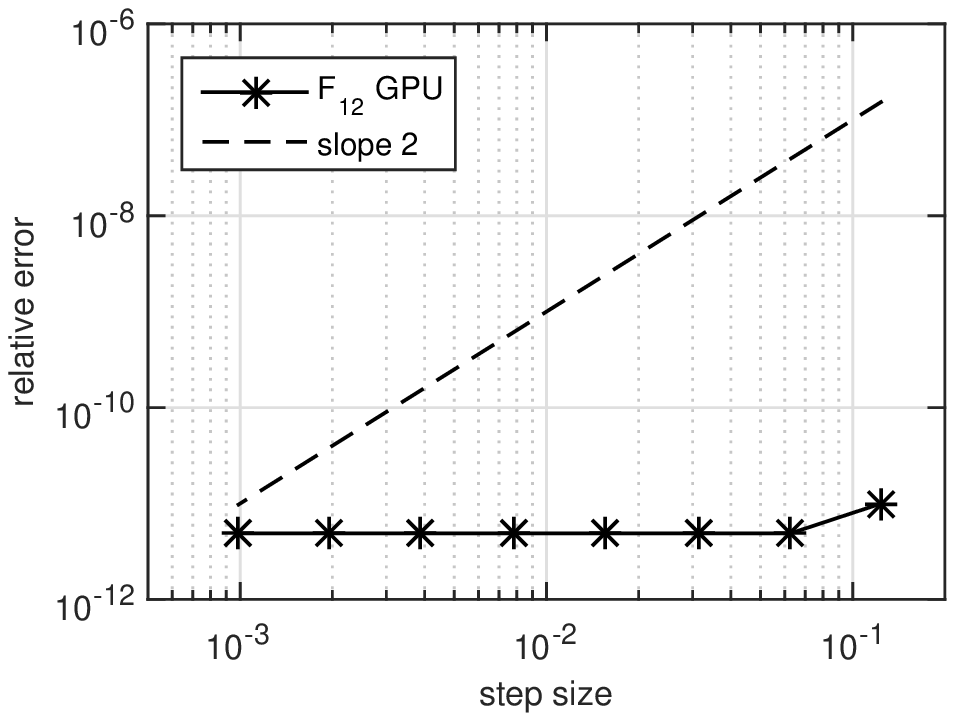}
  \caption{Relative errors of the Strang splitting scheme (left) and the quadrature rule method (right) applied to the DLE in \nameref{example:one}. }
  \label{fig:err1}
\end{figure}

The left plot of Figure~\ref{fig:err1} shows an order plot for the Strang splitting (Algorithm~\ref{alg:DLE_splitting}) applied to the DLE~\eqref{eq:lyap} arising from \nameref{example:one}, and the right plot shows the corresponding results for the quadrature rule method (Algorithm~\ref{alg:DLE_quad_rule}). Here, and in the following, we identify the methods in the figure legends by in which order the subproblems are solved. Thus the splitting in this example is written as $F_1 F_2$ and the quadrature scheme is denoted by $F_{12}$.

We note that the Strang splitting achieves 2nd-order convergence as expected. The quadrature rule, on the other hand, yields a constant but very low error. This is in fact also the expected result, and the error is the error of the quadrature approximation to the integral.
While Strang splitting has recently been suggested multiple times for DLEs, the extra cost for the quadrature rule is in our implementation only $14$ additional evaluations of the matrix exponential action, and we therefore expect that the quadrature rule will essentially always outperform the splitting. This is confirmed in the next section.

\begin{figure}[!htb]
  \includegraphics[width=0.49\textwidth]{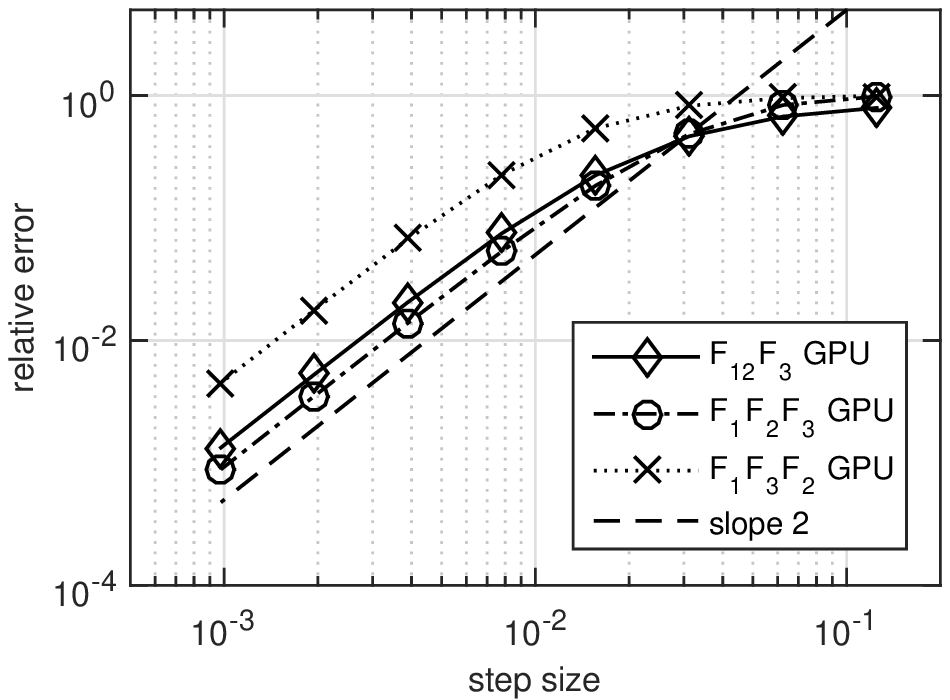}
  \includegraphics[width=0.49\textwidth]{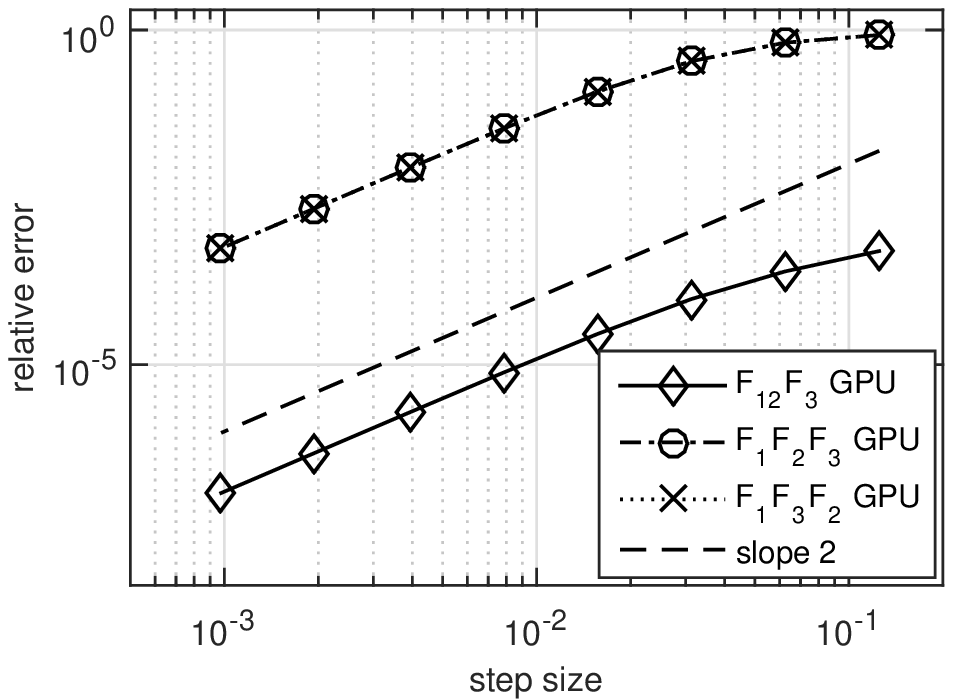}
  \caption{Relative error of the different splitting schemes applied to DRE with $\Ri = 1$ (left) and $\Ri = 10^{-3}$ (right). }
  \label{fig:err_ricc}
\end{figure}

The situation is different for DREs, where we can split into either two or three terms, and the results depend on how large the nonlinear term is compared to the constant term. We consider the three approaches for DREs outlined in Section~\ref{sec:implementation} and denote them by $F_{12}F_3$ (quadrature for the DLE part), $F_1F_2F_3$ and $F_1F_3F_2$. In Figure~\ref{fig:err_ricc} we show an order plot for these methods applied to the DRE~\eqref{eq:ricc} arising from \nameref{example:one}. The left plot uses $\Ri = 1$ and the right one $\Ri = 10^{-3}$.

The first observation to be made is that all the methods converge with the correct order. 
We also see that the three-term splitting $F_1 F_3 F_2$ is less accurate when $\Ri = 1$, whereas the errors of the two remaining splitting schemes behave similarly. Thus, the error due to splitting off the part $F_3$ is more severe than splitting $F_1$ and $F_2$. This is because the nonlinear term is the dominant part here. Using instead $\Ri = 10^{-3}$ means that it is less significant, and leads to a different result. We see that the three-term splittings now yield roughly equally large errors, but that the two-term splitting is about $10$ times more accurate than the other schemes. Here we clearly observe the additional error introduced by the third splitting term.

\begin{figure}[!htb]
	\centering
	\includegraphics[width=0.49\textwidth]{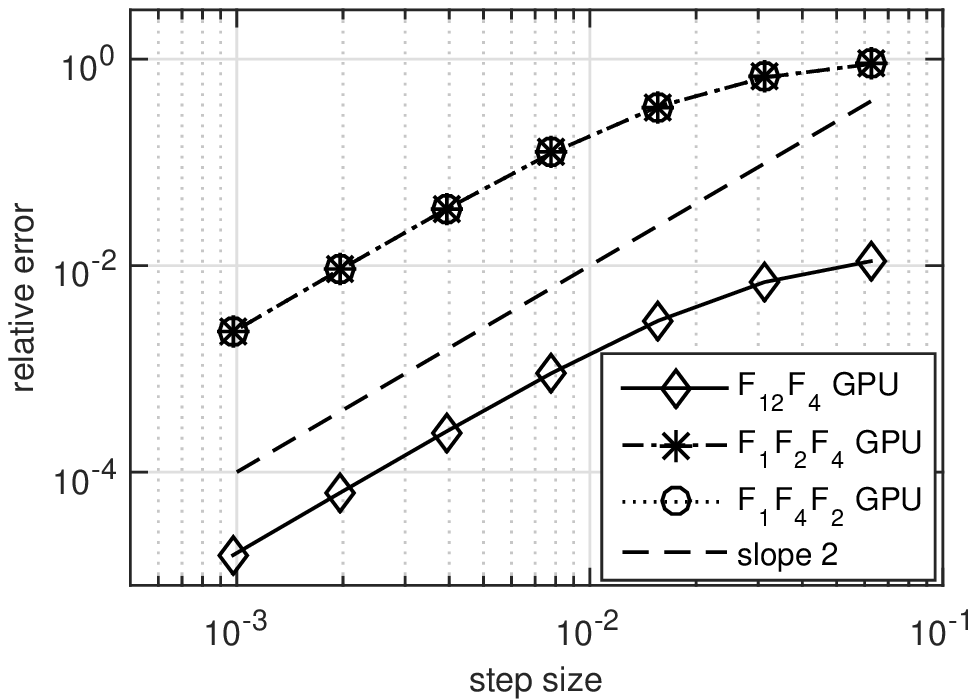}
	\includegraphics[width=0.49\textwidth]{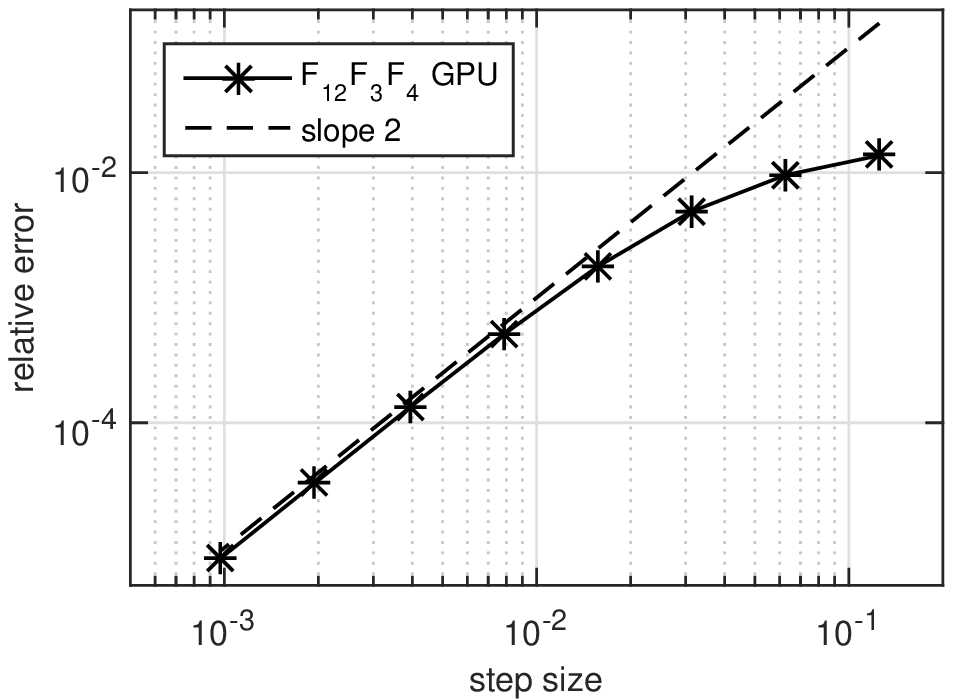}
	\caption{Relative errors of the different splitting schemes applied to the generalized DLE (left) and the generalized DRE (right). }
	\label{fig:err1_gen_eq}
\end{figure}

We also solve the generalized DLE~\eqref{eq:lyapGen} arising from \nameref{example:two} by the three methods mentioned in Section~\ref{sec:implementation} and show the corresponding errors in Figure~\ref{fig:err1_gen_eq} (left). Moreover, we take $R = 1$ and solve also the generalized DRE with the three-term Strang splitting and present the resulting errors in Figure~\ref{fig:err1_gen_eq} (right)
We again see that the two-term splitting of the generalized DLE is approximately 10 times more accurate than the other two splitting schemes. As in all previous examples, we also observe that the error of the generalized DRE behaves as expected, i.e.\ it converges with order 2 and remains small for all step sizes. 

\begin{figure}[!htb]
        \includegraphics[width=0.49\textwidth]{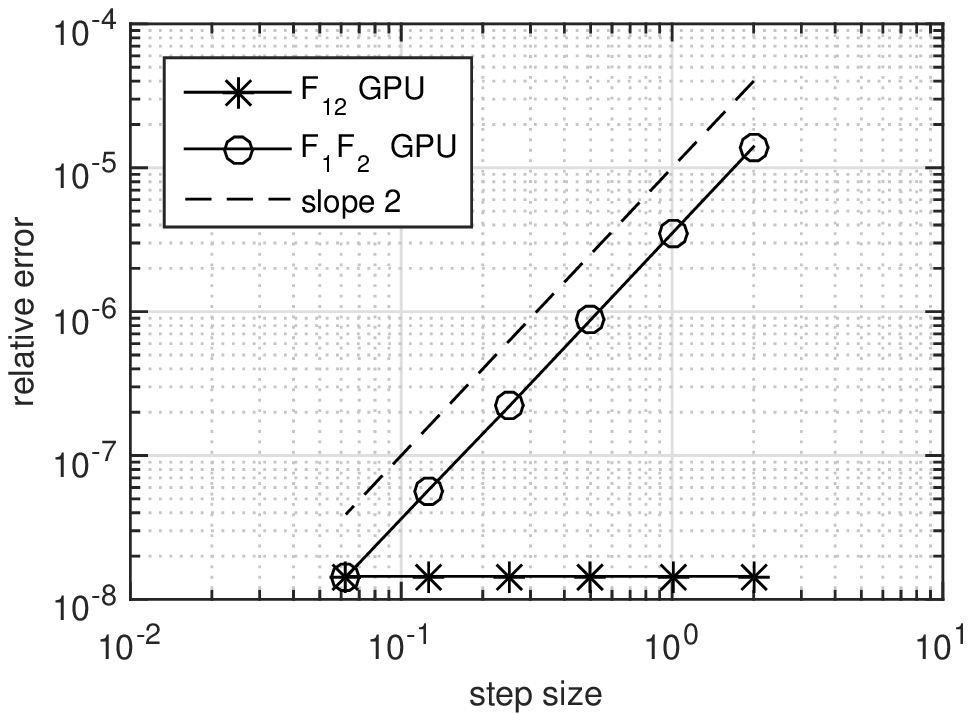}
        \includegraphics[width=0.49\textwidth]{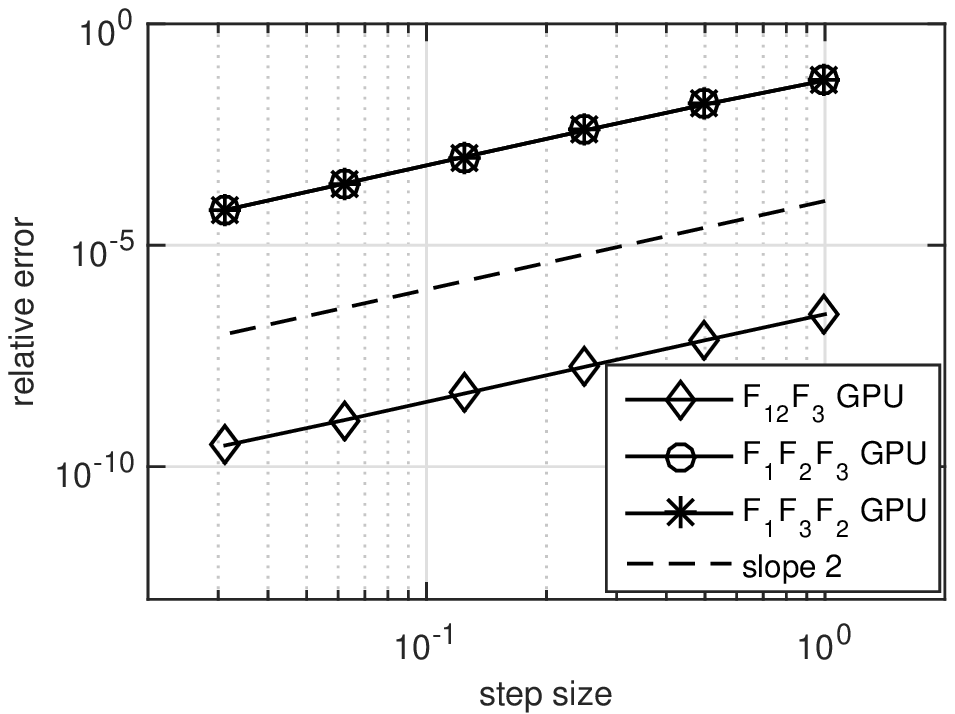}
        \caption{The relative error for the quadrature rule applied to the El Ni\~no DLE with $n =  1740$  (left) and the relative error for three different splitting schemes applied to the steel cooling DRE with $n =  1357$ (right)}
	\label{fig:real_world_order}
\end{figure}

Finally, we test our implementation also on the larger Examples~\ref{example:three} and~\ref{example:four}. Unlike the previous small-scale examples, we do not have a similarly exact reference solution. We instead compare our approximations to an approximation computed with the same scheme, but with a $16$ times smaller step size.
The left plot of Figure~\ref{fig:real_world_order} shows the relative error of the quadrature scheme versus the step sizes when applied to the DLE arising in \nameref{example:three}, and the right plot shows the errors of the two-term and two different three-term splitting schemes applied to the DRE arising from \nameref{example:four}. The problem sizes are here $n = 1740$ and $n = 1357$, respectively, but the error behaves similarly for the other problem sizes. We see that the quadrature rule again produces an essentially exact solution regardless of step size, while the splitting schemes all converge with order two. The two-term splitting once again yields a much lower error than the three-term versions.

\subsection{Performance on main sub-functions}\label{subsec:performance_subfunctions} 

In the following we show the performance of the main sub-functions for the quadrature rule applied to the DLE arising from \nameref{example:one} with $n = 22500$. The other methods and problem cases behave similarly. Using the Matlab profiler, one can show that in the CPU implementation around $98\%$ of the total time is spent in the \texttt{expleja}, which in turn spends almost all its time in the Newton interpolation function. For the GPU implementation the corresponding value is about $90\%$ of the total cost. Hence, we focus on the main sub-functions in this Newton algorithm. They consist of the multiplication of a sparse matrix with a dense skinny matrix (denoted SpMM in the following), the computation of the 1-norm of a dense skinny matrix (1-norm) and the addition of two skinny matrices (addition). The time spent in these main sub-functions, relative to the total computation time, is shown in Table~\ref{table:tab:subfunctions} for the CPU and GPU versions of the code and both systems.
   
\begin{figure}[!htb]
	\includegraphics[width=0.325\textwidth]{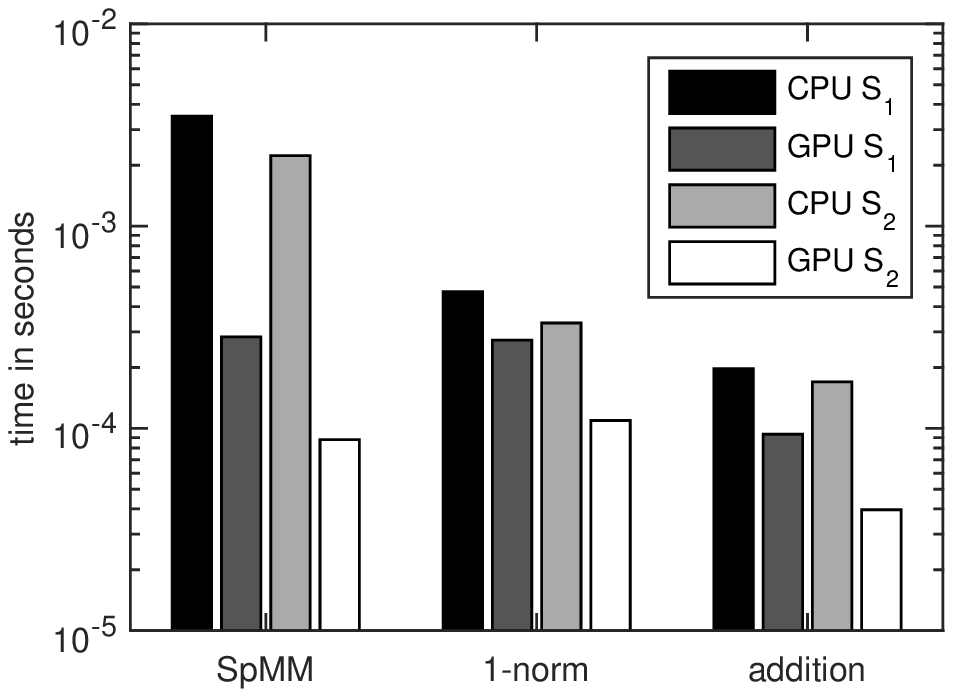}
	\includegraphics[width=0.325\textwidth]{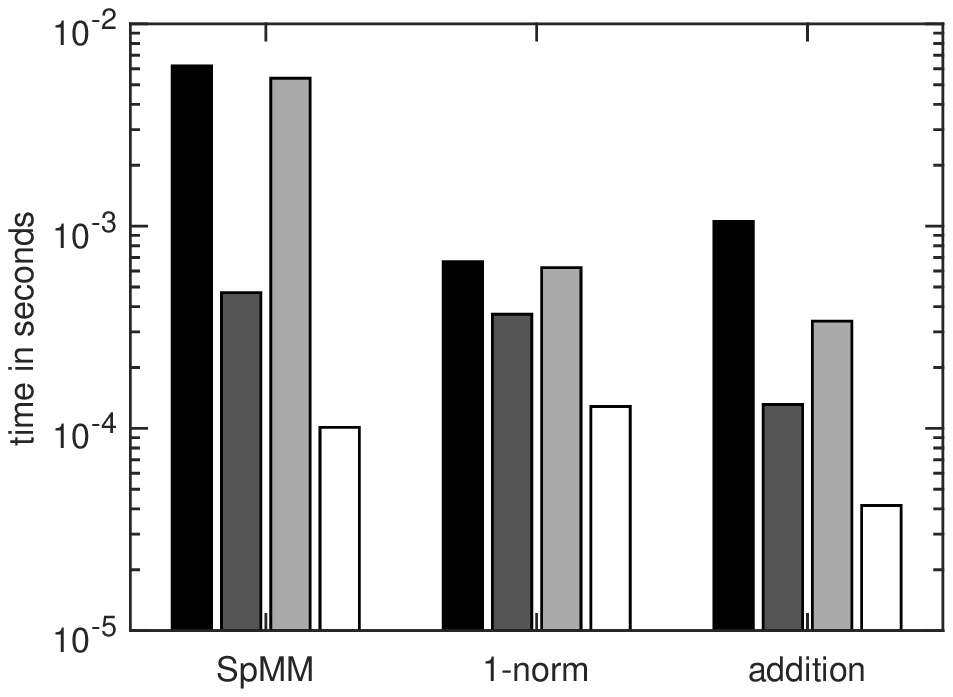}
	\includegraphics[width=0.325\textwidth]{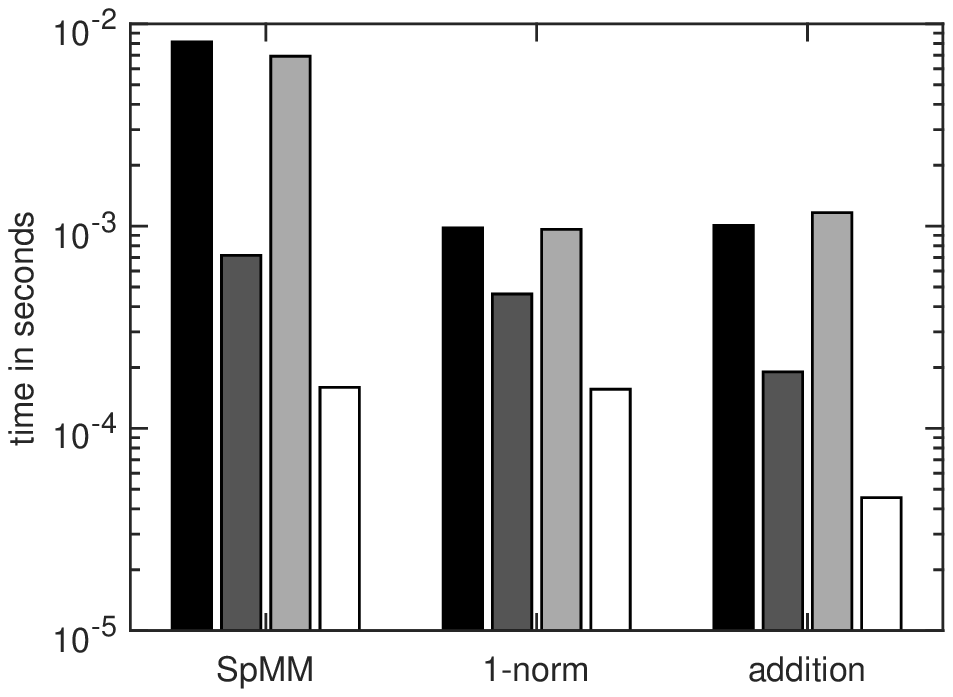}
	\caption{Logarithmic plot of the computational costs of the 3 main operations for matrices with rank 15 (left), rank 30 (middle) and rank 45 (right). }
	\label{fig:main_operations}
\end{figure}

\begin{table}[!htb]                                                           
	\centering                                                                    
	\begin{tabular}{cc|c|c|c|c}
          \multicolumn{2}{c|}{Machine} & Newton & SpMM & 1-norm & addition  \\ 
          \hline\rule{0pt}{3ex}
          \multirow{2}{*}{System 1} & CPU & $98.6\%$ & $72.4\%$ & $19.3\%$   & $4.4\%$ \\   
                                    & GPU & $93.3\%$ & $40.4\%$ & $36.2\%$   & $15.4\%$ \\
          \hline\rule{0pt}{3ex}
          \multirow{2}{*}{System 2} & CPU & $98.3\%$ & $65.2\%$ & $27.1\%$   & $3.8\%$ \\   
                                    & GPU & $89.4\%$ & $35.8\%$ & $36.9\%$   & $14.0\%$ \\   
 	\end{tabular}                                                                 
	\caption{Relative costs of the main sub-functions, in terms of the total computation time.}               
	\label{table:tab:subfunctions}                                                         
\end{table}    
As already mentioned in Subsection~\ref{subsec:possible_gains}, we do not know which algorithm Matlab uses for these sub-functions, but we can compare their costs on the CPU and GPU. We show in Figure~\ref{fig:main_operations} the computational costs for the 3 main sub-functions, applied to randomly generated skinny matrices of ranks $15$, $30$ and $45$. These ranks correspond to the typical ranks of the solutions to the matrix equations arising from \nameref{example:one} and \nameref{example:two}. Here, and in the following, we denote the different systems by $S_1$ and $S_2$ in the figure legends, to save space.

We see that for System 1 we obtain a GPU speed-up of a factor $11-13$ depending on the rank. This is higher than the expected theoretical factor $4.07$. For the computation of the norm, however, we get only a factor of $1.8$. The speed-up of the addition varies between $2.1$ and $8.0$. These different numbers reflect different optimization strategies and uses of cache memory on the different platforms. Since the SpMM multiplication takes up more than $70\%$ of the total costs for the CPU, we draw the conclusion that we can expect a speed-up which is higher than $4$.

For System 2, the SpMM speed-up depends highly on the rank; we get a factor $25.3$ for rank $15$, $ 53.2$ for rank $30$ and $43.4$ for rank $45$. This is again higher than the theoretical factor $6.37$. The speed-up of the 1-norm is between $3.0$ and $6.2$ and the addition varies between $4.2$ and  $25.7$. We thus again expect to see a speed-up higher than what would be expected if both the underlying libraries operated at peak efficicency.

\subsection{Overall performance}
Next, we measure the computational times for the full algorithms. First we consider \nameref{example:one} with the different sizes $n = 625, 2500, 5625, 10000, 15625$ and $22500$, and the step size $h = 0.005$.

\begin{figure}[!htb]
  \includegraphics[width=0.49\textwidth]{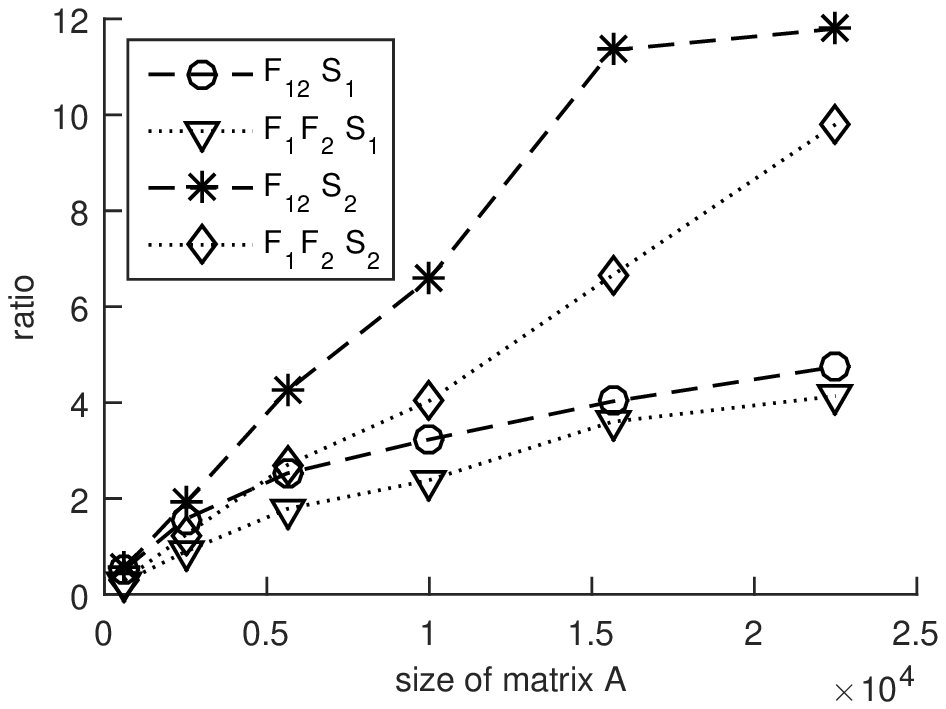}
  \includegraphics[width=0.49\textwidth]{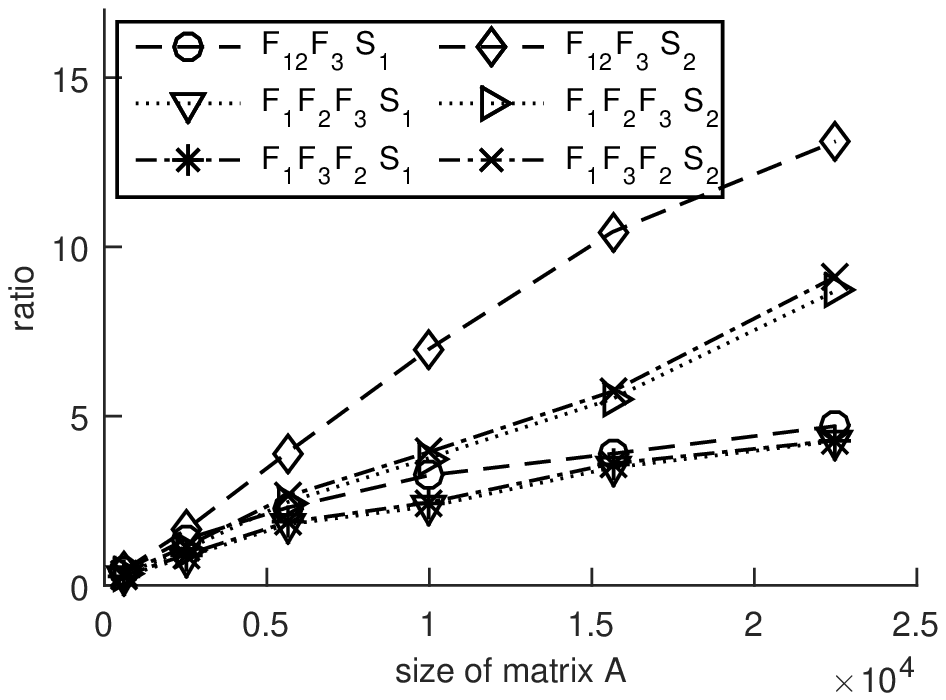}
  \caption{Relative computational costs of the algorithms applied to the DLE (left) and DRE (right) arising from \nameref{example:one}, for different problem sizes. }
  \label{fig:cost_DLEDRE}
\end{figure} 
\begin{table}[!htb]                                                           
	\centering                                                                    
	\begin{tabular}{lr|*{5}{S}}                                                                          
          & $n$ & {2500} & {5625} & {10000} & {15625} & {22500} \\ 
          \hline\rule{0pt}{3ex}
          \multirow{4}{*}{System 1} & $F_{12}$ CPU & 8.310550E+00 & 3.470490E+01 & 8.334515E+01 & 1.924065E+02 & 4.063470E+02 \\   
          & $F_{12}$ GPU & 5.284550E+00 & 1.374630E+01 & 2.579300E+01 & 4.783015E+01 & 8.557915E+01 \\   
          & $F_1 F_2$ CPU & 4.357300E+00 & 2.032330E+01 & 4.623775E+01 & 1.034593E+02 & 1.998573E+02 \\
          & $F_1 F_2$ GPU & 4.877400E+00 & 1.139720E+01 & 1.939640E+01 & 2.874435E+01 & 4.831265E+01 \\
          \hline\rule{0pt}{3ex}
          \multirow{4}{*}{System 2} & 
          $F_{12}$ CPU & 7.562565E+00 & 3.847883E+01 & 9.792389E+01 & 2.672482E+02 & 4.541512E+02 \\   
          & $F_{12}$ GPU & 3.960291E+00 & 9.064024E+00 & 1.487654E+01 & 2.354802E+01 & 3.850605E+01 \\   
          & $F_1 F_2$ CPU & 3.890690E+00 & 1.975971E+01 & 4.800419E+01 & 1.186730E+02 & 2.644458E+02 \\
          & $F_1 F_2$ GPU & 3.136679E+00 & 7.329920E+00 & 1.188841E+01 & 1.784606E+01 & 2.697951E+01 \\
\end{tabular}                                                                                 
\caption{Computational costs of the algorithms applied to the DLE arising from \nameref{example:one}, for different problem sizes.}                               
\label{table:tab:dle}                                                                         
\end{table}

\begin{table}[!htb]                                                           
	\centering                                                                    
	\begin{tabular}{lr|*{5}{S}}                                                                          
          & $n$ & {2500} & {5625} & {10000} & {15625} & {22500} \\ 
          \hline\rule{0pt}{3ex}
          \multirow{6}{*}{System 1} & 
          $F_{12}F_3$ CPU & 8.478950E+00 & 3.534700E+01 & 8.798855E+01 & 1.948567E+02 & 4.125744E+02 \\  
          & $F_{12}F_3$ GPU & 6.214200E+00 & 1.539775E+01 & 2.703465E+01 & 5.009345E+01 & 8.765510E+01 \\  
          & $F_1 F_2 F_3$ CPU & 4.490200E+00 & 2.082990E+01 & 4.630190E+01 & 1.065281E+02 & 2.123421E+02 \\
          & $F_1 F_2 F_3$ GPU & 5.395600E+00 & 1.155230E+01 & 1.984600E+01 & 3.063440E+01 & 4.985260E+01 \\
          & $F_1 F_3 F_2$ CPU & 4.446400E+00 & 2.082280E+01 & 4.623040E+01 & 1.060745E+02 & 2.105331E+02 \\
          & $F_1 F_3 F_2$ GPU & 4.915000E+00 & 1.115160E+01 & 1.886960E+01 & 2.946480E+01 & 4.916400E+01 \\
          \hline\rule{0pt}{3ex}                    
          \multirow{6}{*}{System 2} &
          $F_{12}F_3$ CPU & 7.812395E+00 & 3.833149E+01 & 1.098079E+02 & 2.674168E+02 & 5.265909E+02 \\  
          & $F_{12}F_3$ GPU & 4.763973E+00 & 9.786869E+00 & 1.571181E+01 & 2.562131E+01 & 4.014514E+01 \\  
          & $F_1 F_2 F_3$ CPU & 4.093812E+00 & 2.025942E+01 & 4.808243E+01 & 1.109288E+02 & 2.492502E+02 \\
          & $F_1 F_2 F_3$ GPU & 3.886297E+00 & 8.246683E+00 & 1.276649E+01 & 2.013757E+01 & 2.865047E+01 \\
          & $F_1 F_3 F_2$ CPU & 3.941071E+00 & 2.008789E+01 & 4.769903E+01 & 1.101864E+02 & 2.501578E+02 \\
          & $F_1 F_3 F_2$ GPU & 3.451162E+00 & 7.630086E+00 & 1.209555E+01 & 1.923648E+01 & 2.743647E+01 \\
\end{tabular}                                                                                     
\caption{Computational costs of the algorithms applied to the DRE arising from \nameref{example:one}, for different problem sizes.}                                   
\label{table:tab:dre}                                                                             
\end{table}

In Figure~\ref{fig:cost_DLEDRE} (left) we plot the ratio between the computing time of the CPU and of the GPU as a function of the problem size when applying the different methods to the arising DLE using both systems. The raw data can also be found in Table~\ref{table:tab:dle}, except for the somewhat uninteresting case $n = 625$ which we omit due to space reasons.
We observe that for small matrices the CPU implementation is less time consuming than the GPU implementation. However, the GPU starts to pay off already for problem sizes around $n = 2500$. For the largest test case $n = 22500$ we observe a speed-up of a factor $4.7$ on System 1 and a factor of $11.7$ on System 2 for the quadrature rule. For the splitting scheme, the speed-up is less but not by much. We note that these ratios are higher than the theoretical numbers which one might expect, but fully in line with the analysis in the previous section. We also remark here that the quadrature method is clearly more efficient than the splitting scheme, since the former method produces a much lower error than the latter while their computational costs are very similar.

A similar behaviour can be seen in Figure~\ref{fig:cost_DLEDRE} (right), where we plot the GPU speed-up of the splitting schemes applied to the arising DRE. The break-even point is again around $n = 2500$, as seen in Table~\ref{table:tab:dre} which presents the raw data. For the largest test case, we again observe a factor $4.7$ speed-up for the GPU implementation of the two-term splitting on System 1. On System 2, the corresponding number is $13.1$.
Again, a speed-up of the implementation on the GPU is detected for these problems. The factors for the three-term splittings are both about $4.2$ on System 1 and about $9$ on System 2, which means that all the methods perform better than what might be expected at first glance, due to differently optimized underlying codebases.

We note that the fact that some of the schemes are faster than the others does not mean that they are more efficient, since their errors are also different. By plotting the errors in Figure~\ref{fig:err_ricc} against the computation times, one can observe that the three-term schemes are most efficient for all error levels when $Q$ is relatively small compared to $R^{-1}$, while the two-term splitting is more efficient otherwise. This also holds in general for other problem sizes.

\begin{figure}[!htb]
  \includegraphics[width=0.49\textwidth]{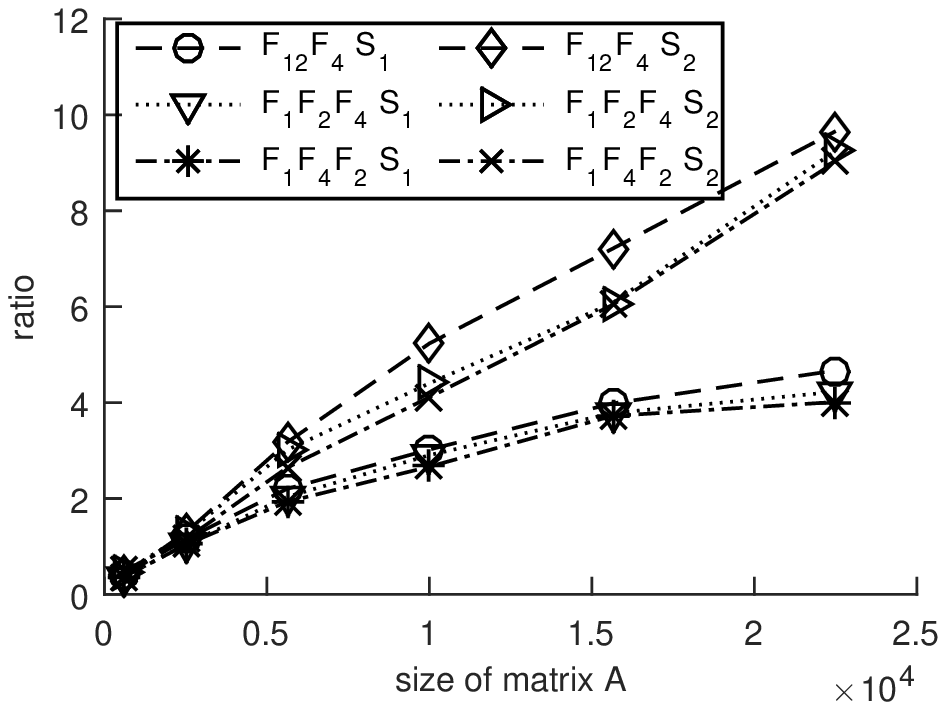}
  \includegraphics[width=0.49\textwidth]{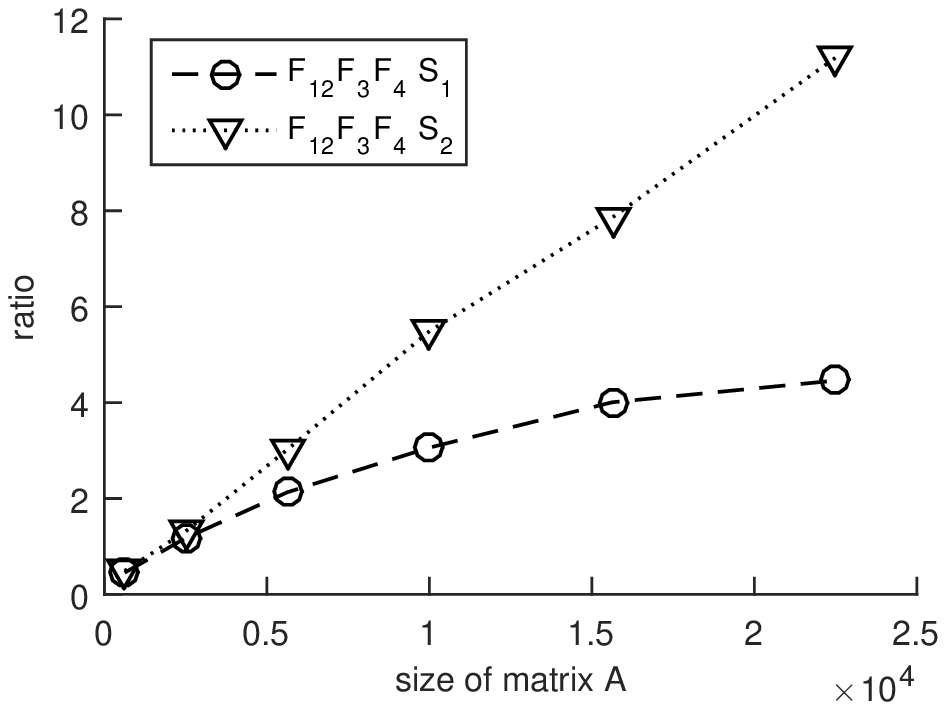}
  \caption{Relative computational costs of the algorithms for the generalized DLE (left) and the generalized DRE (right) arising from  \nameref{example:two}, plotted versus the different problem sizes.}
  \label{fig:cost_gen}
\end{figure}

The GPU-based codes also exhibit better performance for the generalized matrix equations, as seen in Figure~\ref{fig:cost_gen} and Tables~\ref{table:bdle} and~\ref{table:bdre}. We consider here the four schemes mentioned in Section~\ref{sec:splitting_schemes} applied to the generalized DLE and DRE arising from \nameref{example:two}. The break-even point is here slightly lower, but the maximal speed-up factors are similar to the previous examples.

\begin{table}[!htb]                                                           
	\centering                                                                    
	\begin{tabular}{lr|*{5}{S}}                                                                          
          & $n$ & {2500} & {5625} & {10000} & {15625} & {22500} \\ 
          \hline\rule{0pt}{3ex}
          \multirow{6}{*}{System 1} &
            $F_{12}F_4$ CPU & 9.308100E+00 & 3.618990E+01 & 8.331320E+01 & 2.114291E+02 & 4.367187E+02 \\  
          & $F_{12}F_4$ GPU & 8.079000E+00 & 1.643800E+01 & 2.758310E+01 & 5.310190E+01 & 9.354490E+01 \\  
          & $F_1 F_2 F_4$ CPU & 7.683800E+00 & 2.961220E+01 & 6.524490E+01 & 1.634218E+02 & 3.168548E+02 \\
          & $F_1 F_2 F_4$ GPU & 7.358100E+00 & 1.455490E+01 & 2.246390E+01 & 4.323510E+01 & 7.498950E+01 \\
          & $F_1 F_4 F_2$ CPU & 8.455400E+00 & 3.181270E+01 & 7.064090E+01 & 1.837422E+02 & 3.675535E+02 \\
          & $F_1 F_4 F_2$ GPU & 8.166600E+00 & 1.646190E+01 & 2.652470E+01 & 4.937940E+01 & 9.160320E+01 \\
          \hline\rule{0pt}{3ex}
          \multirow{6}{*}{System 2} &
            $F_{12}F_4$ CPU & 9.206597E+00 & 4.029501E+01 & 1.012586E+02 & 2.298941E+02 & 5.017397E+02 \\  
          & $F_{12}F_4$ GPU & 7.150006E+00 & 1.270669E+01 & 1.935123E+01 & 3.190899E+01 & 5.194412E+01 \\  
          & $F_1 F_2 F_4$ CPU & 6.747214E+00 & 3.042367E+01 & 7.052787E+01 & 1.771107E+02 & 3.628452E+02 \\
          & $F_1 F_2 F_4$ GPU & 5.390702E+00 & 1.010383E+01 & 1.602253E+01 & 2.916726E+01 & 3.923785E+01 \\
          & $F_1 F_4 F_2$ CPU & 7.393654E+00 & 3.296680E+01 & 8.333982E+01 & 1.934752E+02 & 4.363592E+02 \\
          & $F_1 F_4 F_2$ GPU & 6.433725E+00 & 1.246483E+01 & 2.026411E+01 & 3.206172E+01 & 4.837870E+01 \\
\end{tabular}                                                                                     
\caption{Computational costs of the algorithms applied to the generalized DLE arising from \nameref{example:two}, for different problem sizes.}
\label{table:bdle}                                                                            
\end{table}

\begin{table}[!htb]                                                                                
	\centering                                                                                         
	\begin{tabular}{lr|*{5}{S}}                                                                          
          & $n$ & {2500} & {5625} & {10000} & {15625} & {22500} \\ 
          \hline\rule{0pt}{3ex}                           
          \multirow{2}{*}{System 1} &
            $F_{12}F_3F_4$ CPU & 9.419900E+00 & 3.652780E+01 & 8.784420E+01 & 2.186546E+02 & 4.331563E+02 \\
          & $F_{12}F_3F_4$ GPU & 8.018000E+00 & 1.715720E+01 & 2.869800E+01 & 5.458100E+01 & 9.713390E+01 \\
          \hline\rule{0pt}{3ex}                           
          \multirow{2}{*}{System 2} &
            $F_{12}F_3F_4$ CPU & 8.721255E+00 & 3.927012E+01 & 1.065744E+02 & 2.496264E+02 & 5.508864E+02 \\
          & $F_{12}F_3F_4$ GPU & 6.629428E+00 & 1.301385E+01 & 1.942490E+01 & 3.177167E+01 & 4.921916E+01 \\
	\end{tabular}                                                                                      
	\caption{Computational costs of the algorithms applied to the generalized DRE arising from \nameref{example:two}, for different problem sizes.}
	\label{table:bdre}                                                                             
\end{table} 

\begin{figure}[!htb]
  \includegraphics[width=0.49\textwidth]{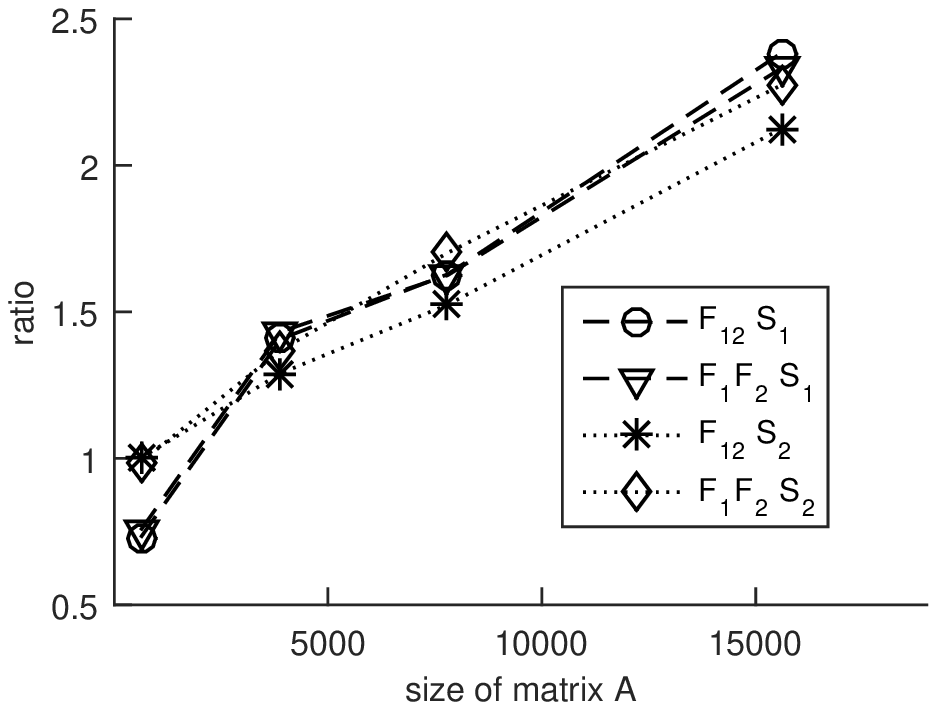}
  \includegraphics[width=0.49\textwidth]{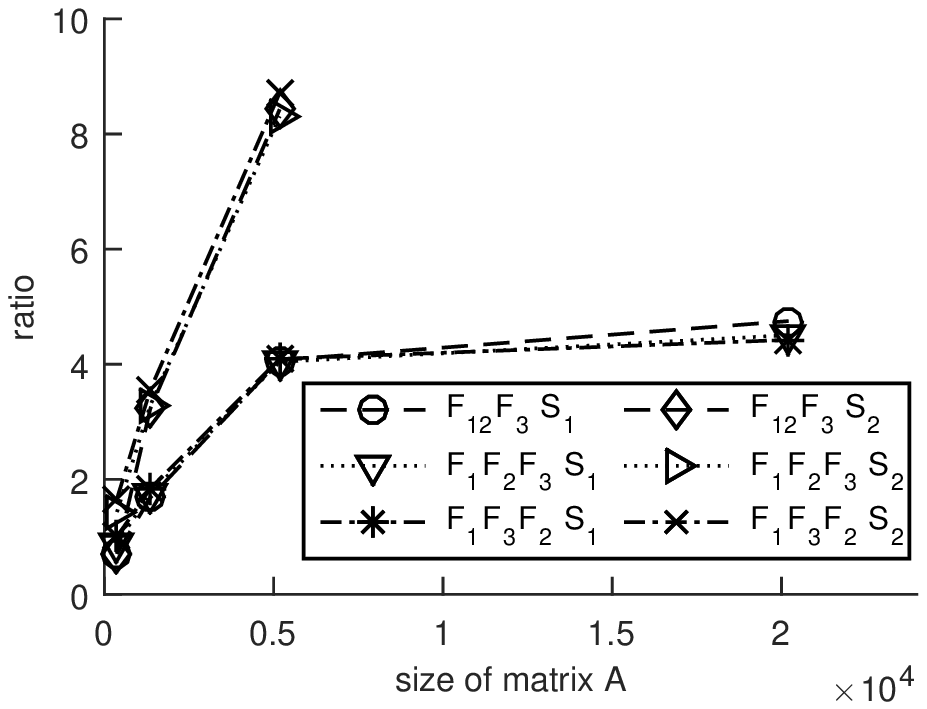}
  \caption{Relative computational costs for \nameref{example:three} (left) and \nameref{example:four} (right). Due to time constraints, we could not test the largest problem size in \nameref{example:four} on System~2.}
  \label{fig:real_world_times}
\end{figure}

Finally, we measure the computation times also for the two real-world examples. The left plot in Figure~\ref{fig:real_world_times} shows the results of applying the DLE methods to \nameref{example:three} for different problem sizes, and the right plot shows the DRE methods applied to \nameref{example:four}. Tables~\ref{table:nino} and~\ref{table:steel} contains the respective raw data.

\begin{table}[!htb]                                                            
	\centering     
        \begin{tabular}{lr|*{4}{S}}                                                                          
          & $n$ & {624} & {3900} & {7800} & {15600} \\               
          \hline\rule{0pt}{3ex}
          \multirow{4}{*}{System 1} &
            $F_{12}$ CPU & 2.552500E+00 & 1.027470E+01 & 1.871200E+01 & 5.020190E+01 \\   
          & $F_{12}$ GPU & 3.498100E+00 & 7.300100E+00 & 1.149020E+01 & 2.106380E+01 \\   
          & $F_1 F_2$ CPU & 2.473100E+00 & 9.821500E+00 & 1.764520E+01 & 4.760840E+01 \\
          & $F_1 F_2$ GPU & 3.277200E+00 & 6.857500E+00 & 1.085760E+01 & 2.042000E+01 \\
          \hline\rule{0pt}{3ex}                                                                  
          \multirow{4}{*}{System 2} &
            $F_{12}$ CPU & 4.086157E+00 & 1.147969E+01 & 2.035560E+01 & 4.631066E+01 \\   
          & $F_{12}$ GPU & 4.081647E+00 & 8.894490E+00 & 1.335339E+01 & 2.179319E+01 \\   
          & $F_1 F_2$ CPU & 5.428569E+00 & 1.093816E+01 & 1.948929E+01 & 4.670509E+01 \\
          & $F_1 F_2$ GPU & 5.501126E+00 & 7.987852E+00 & 1.145281E+01 & 2.053894E+01 \\
\end{tabular}                                                                             
	\caption{Computational costs of the algorithms applied to the DLE arising from \nameref{example:three}, for different problem sizes.}                
	\label{table:nino}                                                          
\end{table}

\begin{table}[!htb]                                                                
	\centering                                                                         
	\begin{tabular}{lr|*{4}{S}}                                                                          
          & $n$ & {371} & {1357} & {5177} & {20209} \\
          \hline\rule{0pt}{3ex}
          \multirow{6}{*}{System 1} &
            $F_{12}F_3$ CPU & 3.228270E+01 & 3.107412E+02 & 2.715800E+03 & 9.222164E+04 \\  
          & $F_{12}F_3$ GPU & 4.565400E+01 & 1.846400E+02 & 6.673000E+02 & 1.941753E+04 \\  
          & $F_1 F_2 F_3$ CPU & 3.367350E+01 & 3.104393E+02 & 2.593100E+03 & 7.884768E+04 \\
          & $F_1 F_2 F_3$ GPU & 3.808600E+01 & 1.767700E+02 & 6.416200E+02 & 1.750390E+04 \\
          & $F_1 F_3 F_2$ CPU & 3.382860E+01 & 3.150115E+02 & 2.557600E+03 & 7.480704E+04 \\
          & $F_1 F_3 F_2$ GPU & 3.372000E+01 & 1.711300E+02 & 6.248900E+02 & 1.691200E+04 \\
          \hline\rule{0pt}{3ex}
          \multirow{6}{*}{System 2} &
		$F_{12}F_3$ CPU  & 4.534219E+01 & 3.400065E+02 & 2.873341E+03 &  \\  
		&$F_{12}F_3$ GPU  & 6.431840E+01 & 1.052602E+02 & 3.409001E+02 &  \\  
		&$F_1 F_2 F_3$ CPU  & 3.656979E+01 & 3.455842E+02 & 2.753539E+03 &  \\
		&$F_1 F_2 F_3$ GPU  & 2.552956E+01 & 1.051054E+02 & 3.315365E+02 &  \\
		&$F_1 F_3 F_2$ CPU  & 3.510726E+01 & 3.419044E+02 & 2.749892E+03 &  \\
		&$F_1 F_3 F_2$ GPU  & 2.113769E+01 & 9.604626E+01 & 3.155070E+02 &  \\
	\end{tabular}                                                                      
	\caption{Computational costs of the algorithms applied to the DRE arising from \nameref{example:four}, for different problem sizes. Due to time constraints, we could not test the largest problem size on System~2.}                    
	\label{table:steel}                                                            
\end{table}                     

The results are similar to the previous academic examples. In the DLE case, the CPU and GPU costs are comparable at $n = 624$, but at $n = 3900$ the GPU is more efficient and at $n = 7800$ we already observe a speed-up of a factor $1.5$. For the finest resolution, the speed-up is roughly a factor $2.4$. The maximal speed-up is lower in this example, because the solution is of a higher rank than in the academic examples. This requires more work in the column compression step, which in turn performs SVD calculations. These are harder to parallelize than the other main sub-functions. We note, however, that as the problem dimension increases the solution rank increases only marginally. This means that for large enough problems the column compression cost will again be negligible, and the GPU speed-up will reach similar values as in the academic examples. We still want to emphasize that for the current largest test case the algorithm performs twice as good on the GPU as on the CPU, and this is with essentially no changes to the code.

In the DRE case, only the first problem size yields comparable costs for the CPU and GPU, and we observe speed-ups for all larger problem sizes. We obtain a speed-up of $4$ already for $n = 5177$ and $4.7$ for the largest test case on System 1. Due to time constraints, we only solve the first three problem sizes on System 2. For $n = 5177$ the speed-up is already more than $8.3$ for all schemes, and we expect the ratio to level out similarly to what happens on System 1. As mentioned previously, even higher speed-up are expected when Matlab supports solving equation systems with sparse system matrices and dense block right-hand sides.

\section{Conclusions} \label{sec:conclusions}
We have considered several different splitting schemes based on Leja point interpolation for the computation of matrix exponential actions. Since the matrix exponentials act on skinny block-matrices (the low-rank factors) rather than only vectors, we expected that these computations would be highly parallelizable and that GPU acceleration would therefore be beneficial. The latter was verified by several numerical experiments on two different systems. In the considered problems of academical nature, the GPU code was faster than the pure CPU code by approximately a factor $3$ already for matrices of size $10000$ on System 1 and by a factor of $6$ on System 2. This factor increases to over $4$ and $10$, respectively, for larger matrices of size $22500$. The break-even point was around size $2500$, which is well below what would be considered large-scale today.In the tested real-world applications, the gains were also in accordance with the more academic examples. As there is no difference in the size of the numerical errors, this clearly shows that GPU acceleration can lead to large gains in efficiency and should be considered for matrix equations of these types. The efficiency could additionally be further increased by considering more advanced parallelization techniques. An obvious such candidate is to investigate the use of single-precision computations when the desired level of accuracy is low.

We have also presented comparisons of different splitting strategies, mainly investigating whether one should split off the constant term $Q$ or not, and in which order the subproblems should be solved. For the latter question, we observe that the ordering has minimal influence on the error, and we may thus choose the order such that the computational cost is minimized. (E.g., take the most expensive subproblem as the ``middle'' term.) For the first question, we expected that it would not be beneficial to split off $Q$, since the extra integral term which arises only has to be approximated once. This was verified by our experiments, except in the case when $Q$ was relatively small -- then, of course, the extra splitting error is similarly small. We note that these results are for the autonomous case. When the matrices that define the equations also depend on time, the situation likely changes, as the integral term would need to be recomputed in each step. However, as the modified methods would still rely on matrix exponential actions as their basic building blocks, we still expect that GPU acceleration would significantly increase the efficiency.

\section*{\ackname}
The authors are grateful to Peter Kandolf for his assistance with the original \texttt{expleja} code. H. Mena and L.-M. Pfurtscheller were supported by the Austrian Science Fund (FWF) - project id:P27926.
%\end{ackname}

\bibliography{thebib}
\bibliographystyle{amsplain}

\end{document}